# On Itô's formula for elliptic diffusion processes

XAVIER BARDINA[1] and CARLES ROVIRA[2]

[1]*Departament de Matemàtiques, Universitat Autònoma de Barcelona, 08193-Bellaterra (Barcelona), Spain. E-mail: Xavier.Bardina@uab.cat*

[2]*Facultat de Matemàtiques, Universitat de Barcelona, Gran Via 585, 08007-Barcelona, Spain. E-mail: carles.rovira@ub.edu*

Bardina and Jolis [*Stochastic process. Appl.* **69** (1997) 83–109] prove an extension of Itô's formula for $F(X_t,t)$, where $F(x,t)$ has a locally square-integrable derivative in $x$ that satisfies a mild continuity condition in $t$ and $X$ is a one-dimensional diffusion process such that the law of $X_t$ has a density satisfying certain properties. This formula was expressed using quadratic covariation. Following the ideas of Eisenbaum [*Potential Anal.* **13** (2000) 303–328] concerning Brownian motion, we show that one can re-express this formula using integration over space and time with respect to local times in place of quadratic covariation. We also show that when the function $F$ has a locally integrable derivative in $t$, we can avoid the mild continuity condition in $t$ for the derivative of $F$ in $x$.

## 1. Introduction

Itô's formula is one of the most important results in the theory of stochastic processes because it plays the role of the theorem of change of variable for stochastic integration.

Itô [11] obtained this result in 1944 for a standard Brownian motion $W = (W_t)_{t \geq 0}$ and a $\mathcal{C}^{2,1}$ function $F: \mathbb{R} \times \mathbb{R}_+ \longrightarrow \mathbb{R}$. Precisely, he proved that

$$F(W_t,t) = F(W_0,0) + \int_0^t \frac{\partial F}{\partial t}(W_s,s)\,\mathrm{d}s + \int_0^t \frac{\partial F}{\partial x}(W_s,s)\,\mathrm{d}W_s + \frac{1}{2}\int_0^t \frac{\partial^2 F}{\partial x^2}(W_s,s)\,\mathrm{d}s.$$

Since then, many extensions of this result have been obtained, in two directions: proving that the formula is also valid for other types of processes and weakening the conditions on the function $F$.

In the first direction, Kunita and Watanabe [12] proved the result for continuous semimartingales $X$ and functions $F \in \mathcal{C}^{2,1}$. The first appearance of local time in an Itô formula is due to Tanaka [15] in 1963, but in 1981, Bouleau and Yor [3] obtained an extension of Itô's formula where the additional term turns up as an integral in space with respect to the local time.

In 1995, Föllmer *et al.* [9] made one of the more remarkable contributions to the theory in the sense of weakening the conditions on the function $F$. Their Itô formula is for $F(W_t, t)$, where $F(x, t)$ has a locally square-integrable derivative in $x$ that satisfies a mild continuity condition in $t$, and $W$ is a standard Brownian motion.

Eisenbaum [5, 6] defined an integral in time and space with respect to the local time of the Brownian motion. Using this integral, the quadratic covariation in the formula given in Föllmer *et al.* can be expressed as an integral with respect to the local time.

Bardina and Jolis [2] extended the results of Föllmer *et al.* [9] to the case of elliptic diffusions. In this paper, our goal is to generalize the results obtained by Eisenbaum to this case.

Other extensions of Itô's formula include Errami *et al.* [8], Dupoiron *et al.* [4], Peskir [14], Ghomrasni and Peskir [10] and Eisenbaum [7].

The paper is organized as follows. In Section 2, we recall the basic results obtained in Bardina and Jolis [2], where the authors obtain an Itô formula for elliptic diffusion processes. In Section 3, we define the space where we are able to construct an integral in the plane with respect to the local time of a diffusion process. Section 4 is devoted to presenting the extension of Itô's formula. Finally, in Section 5, we recall some examples of diffusions satisfying the hypotheses considered in Section 4.

## 2. Preliminaries

Bardina and Jolis [2] prove an extension of Itô's formula for $F(X_t, t)$, where $F(x, t)$ has a locally square-integrable derivative in $x$ that satisfies a mild continuity condition in $t$, and $X$ is a one-dimensional diffusion process such that the law of $X_t$ has a density satisfying certain properties.

Following the ideas of Föllmer *et al.* [9], where they prove an analogous extension when $X$ is Brownian motion, the proof is based on the existence of the quadratic covariation $[f(X, \cdot), X]_t$ defined as the limit uniformly in probability of

$$\sum_{t_i \in D_n, t_i \leq t} (f(X_{t_{i+1}}, t_{i+1}) - f(X_{t_i}, t_i))(X_{t_{i+1}} - X_{t_i}),$$

where $(D_n)_n$ denotes a sequence of partitions of $[0, 1]$, each composed of intervals with the norm tending to zero and all having the same order of convergence. More precisely, they prove that the limit uniformly in probability of

$$\sum_{t_i \in D_n, t_i \leq t} f(X_{t_i}, t_i)(X_{t_{i+1}} - X_{t_i})$$

exists and coincides with the forward stochastic integral $\int_0^t f(X_s, s) \, dX_s$ and that the limit uniformly in probability of

$$\sum_{t_i \in D_n, t_i \leq t} f(X_{t_{i+1}}, t_{i+1})(X_{t_{i+1}} - X_{t_i})$$

exists and coincides with the backward stochastic integral $\int_0^t f(X_s, s) \, \mathrm{d}^* X_s$. The backward integral is defined as a forward stochastic integral of $f(\bar{X}_s, 1-s)$ with respect to the time-reversed process $\bar{X} = \{\bar{X}_t \equiv X_{1-t}, 0 \leq t \leq 1\}$:

$$\int_0^t f(X_s, s) \, \mathrm{d}^* X_s = - \int_{1-t}^1 f(\bar{X}_s, 1-s) \, \mathrm{d}\bar{X}_s.$$

The backward integral is well defined when the reversed process satisfies certain integrator properties. For diffusion processes, there are results that guarantee that the reversed process is also a diffusion and therefore a semimartingale.

Assume that our process satisfies the stochastic integral equation

$$X_t = X_0 + \int_0^t b(s, X_s) \, \mathrm{d}s + \int_0^t \sigma(s, X_s) \, \mathrm{d}W_s,$$

where $W = \{W_t, \mathcal{F}_t; 0 \leq t \leq 1\}$ is a standard Brownian motion, and $\sigma$ and $b$ satisfy the usual conditions (H) of Lipschitz and linear growth: there exists a constant $K > 0$ such that for any $x, y \in \mathbb{R}$,

$$\sup_t [|\sigma(t, x) - \sigma(t, y)| + |b(t, x) - b(t, y)|] \leq K|x - y|,$$

$$\sup_t [|\sigma(t, x)| + |b(t, x)|] \leq K(1 + |x|).$$

The reversibility of the diffusion property has been studied by several authors; see, for instance, Millet et al. [13]. The extension of Itô's formula obtained for Bardina and Jolis is as follows.

**Theorem 2.1 (Bardina and Jolis [2]).**

(H1) *Hypotheses on the function:*

1. *Let $F(x, t)$ be an absolutely continuous function with respect to the coordinate $x$ and let $f(\cdot, t)$ be its derivative $\frac{\partial F}{\partial x}(\cdot, t)$.*
2. *Suppose that $f(\cdot, t) \in L^2_{loc}(\mathbb{R})$ and that for all compact subsets $T$ of $\mathbb{R}$, $\int_T f^2(x, t) \, \mathrm{d}x$ is a continuous function of $t$ in $[0, 1]$.*

(H2) *Hypotheses on the diffusion:*

1. *Assume that $X_t$ is a diffusion process that satisfies hypotheses* (H).
2. *Suppose that for all $t > 0$, $X_t$ has a density $p_t(x)$ that satisfies the two conditions*

$$\sup_x p_t(x) \leq r(t),$$

$$\left( \int_{\mathbb{R}} \frac{1}{p_t(x)} \left( \frac{\mathrm{d}}{\mathrm{d}x} (\sigma^2(t, x) p_t(x)) \right)^2 \mathrm{d}x \right)^{1/2} \leq u(t),$$

*where $r(t)$, $u(t)$, $t \in (0, 1]$ are continuous non-increasing integrable functions. (The derivative $\frac{\mathrm{d}}{\mathrm{d}x}(\sigma^2(t, x) p_t(x))$ must be interpreted in the distributional sense and, by convention, $(p_t(x))^{-1} = 0$ if $p_t(x) = 0$.)*

3. Assume that $(r)^{1/2}u \in L^1[0,1]$.

*Then, the quadratic covariation exists, as a limit uniformly in probability, and*

$$[f(X,\cdot),X]_t = \int_0^t f(X_s,s)\,\mathrm{d}^*X_s - \int_0^t f(X_s,s)\,\mathrm{d}X_s.$$

*Moreover,*

$$F(X_t,t) = F(X_0,0) + \int_0^t f(X_s,s)\,\mathrm{d}X_s + \tfrac{1}{2}[f(X,\cdot),X]_t + \int_0^t F(X_s,\mathrm{d}s),$$

*where*

$$\int_0^t F(X_s,\mathrm{d}s) \equiv \lim_{n\to+\infty} \sum_{t_i\in D_n, t_i\leq t} (F(X_{t_{i+1}},t_{i+1}) - F(X_{t_{i+1}},t_i)),$$

*exists uniformly in probability.*

**Remark 2.2.** Notice that in Theorem 2.1, if we assume stronger conditions on the function, we can weaken the hypotheses on the diffusion. Particularly, when $f$ is bounded, the condition $(r)^{1/2}u \in L^1[0,1]$ is not necessary.

**Remark 2.3.** In Theorem 2.1, if $F$ is of class $\mathcal{C}^1$ in $t$, then the term $\int_0^1 F(X_s,\mathrm{d}s)$ equals $\int_0^t \frac{\partial F}{\partial s}(X_s,s)\,\mathrm{d}s$. On the other hand, if $F$ is of class $\mathcal{C}^2$ in $x$, then $[f(X,\cdot),X]_t = \int_0^t \frac{\partial f}{\partial x}(X_s,s)\,\mathrm{d}\langle X,X\rangle_s$.

**Remark 2.4.** Under hypotheses (H2) (see Theorem 2.1) and using Theorem 2.3 of Millet *et al.* [13] (also see Remarks 1.2, 1.3 and Theorem 2.2 of Bardina and Jolis [2]), we have the following semimartingale expression of $\bar{X}$ for $0 \leq t \leq 1$:

$$\bar{X}_t = \bar{X}_0 + \int_0^t \bar{b}(s,\bar{X}_s)\,\mathrm{d}s + \int_0^t \bar{\sigma}(s,\bar{X}_s)\,\mathrm{d}\bar{W}_s,$$

where $\{\bar{W}_t, 0 \leq t \leq 1\}$ is another standard Brownian motion defined in the same probability space as $W$ and where

$$\bar{\sigma}(1-t,x) = \sigma(t,x),$$
$$\bar{b}(1-t,x) = -b(t,x) + \frac{1}{p_t(x)}\frac{\mathrm{d}}{\mathrm{d}x}(\sigma^2(t,x)p_t(x)),$$

with the convention that $(p_t(x))^{-1} = 0$ if $p_t(x) = 0$.

Finally, we will represent by $K$ the constants that will appear, although they may vary from occurrence to occurrence.

## 3. Stochastic integration with respect to local time of the diffusion process

Following the ideas of Eisenbaum [5], we first consider the space of functions for whose elements we can define a stochastic integration with respect to local time.

Let us consider $X$, a diffusion process satisfying the stochastic integral equation

$$X_t = X_0 + \int_0^t b(s, X_s)\,\mathrm{d}s + \int_0^t \sigma(s, X_s)\,\mathrm{d}W_s.$$

Note that we need the coefficients $\sigma$ and $b$ to satisfy the hypotheses in order for the reversed time process to also be a martingale. In fact, we will suppose a little more and we will assume statements 1 and 2 of (H2).

Let $f$ be a measurable function from $\mathbb{R} \times [0,1]$ to $\mathbb{R}$. We define the seminorm $\|\cdot\|_*$ by

$$\|f\|_* = 2\int_0^1 \int_\mathbb{R} |f(x,s)b(s,x)|p_s(x)\,\mathrm{d}x\,\mathrm{d}s + 2\left(\int_0^1 \int_\mathbb{R} f^2(x,s)\sigma^2(s,x)p_s(x)\,\mathrm{d}x\,\mathrm{d}s\right)^{1/2}$$
$$+ \int_0^1 \int_\mathbb{R} \left|f(x,s)\frac{\mathrm{d}}{\mathrm{d}x}(\sigma^2(s,x)p_s(x))\right|\mathrm{d}x\,\mathrm{d}s,$$

where the derivative of $\sigma^2(s,x)p_s(x)$ must be interpreted in the distributional sense.

We also assume some hypotheses on the diffusion in order to ensure that the seminorm $\|\cdot\|_*$ is a norm, that is,

(H3) $p_s(x) \neq 0$ a.s. and $\lambda\{(x,t): \sigma(x,t) = b(x,t) = 0\} = 0$, where $\lambda$ denotes the Lebesgue measure.

Consider the set of functions

$$\mathcal{H} = \{f : \|f\|_* < +\infty\}.$$

It is easy to check that $\mathcal{H}$ is a Banach space.

In Remark 3.3, we will show how to define the norm and the space when (H3) is not satisfied, in order to be able to construct the stochastic integration with respect to the local time.

*Remark 3.1.* Note that for Brownian motion, that is, when $b \equiv 0$ and $\sigma \equiv 1$, we obtain the norm and the space given by Eisenbaum [5].

Let us now show how to define a stochastic integration over the plane with respect to the local time $L$ of the diffusion process $X$ for the elements of $\mathcal{H}$.

Let $f_\Delta$ be an elementary function,

$$f_\Delta(x,s) = \sum_{(x_i, s_j) \in \Delta} f_{ij} I_{(x_i, x_{i+1}]}(x) I_{(s_j, s_{j+1}]}(s),$$

where $(x_i)_{1\leq i\leq n}$ is a finite sequence of real numbers, $(s_j)_{1\leq j\leq m}$ is a subdivision of $[0,1]$, $(f_{ij})_{1\leq i\leq n; 1\leq j\leq m}$ is a sequence of real numbers and, finally, $\Delta = \{(x_i, s_j), 1\leq i\leq n, 1\leq j\leq m\}$. It is easy to check that the elementary functions are dense in $\mathcal{H}$.

We define the integration of the elementary function $f_\Delta$ with respect to the local time $L$ of the diffusion $X$ as follows:

$$\int_0^1 \int_\mathbb{R} f_\Delta(x,s)\,\mathrm{d}L_s^x = \sum_{(x_i,s_j)\in\Delta} f_{ij}(L_{s_{j+1}}^{x_{i+1}} - L_{s_j}^{x_{i+1}} - L_{s_{j+1}}^{x_i} + L_{s_j}^{x_i}).$$

Let $f$ be a function of $\mathcal{H}$. Let us consider $(f_n)_{n\in\mathbb{N}}$, a sequence of elementary functions converging to $f$ in $\mathcal{H}$. We will check that the sequence $(\int_0^1 \int_\mathbb{R} f_n(x,s)\,\mathrm{d}L_s^x)_{n\in\mathbb{N}}$ converges in $L^1$ and that the limit does not depend on the choice of sequence $(f_n)_{n\in\mathbb{N}}$. We will thus use this limit as the definition of the integral $\int_0^1 \int_\mathbb{R} f(x,s)\,\mathrm{d}L_s^x$.

**Theorem 3.2.** *Let $f$ be a function of $\mathcal{H}$. Then the integral $\int_0^t \int_\mathbb{R} f(x,s)\,\mathrm{d}L_s^x$ exists and*

$$\int_0^t \int_\mathbb{R} f(x,s)\,\mathrm{d}L_s^x = \int_0^t f(X_s,s)\,\mathrm{d}X_s - \int_0^t f(X_s,s)\,\mathrm{d}^*X_s,$$

*for any $t \in [0,1]$.*

**Proof.** Let $f_\Delta$ be an elementary function. Following the ideas of the proof of Theorem 3.1 of Eisenbaum [5], it is easy to show that

$$\int_0^t \int_\mathbb{R} f_\Delta(x,s)\,\mathrm{d}L_s^x = \int_0^t f_\Delta(X_s,s)\,\mathrm{d}X_s - \int_0^t f_\Delta(X_s,s)\,\mathrm{d}^*X_s.$$

Using the semimartingale expression of the reversed process $\bar{X}$ given in Remark 2.4, we can write

$$\int_0^t f_\Delta(X_s,s)\,\mathrm{d}^*X_s = \int_{1-t}^1 f_\Delta(\bar{X}_s, 1-s)b(1-s,\bar{X}_s)\,\mathrm{d}s$$
$$- \int_{1-t}^1 f_\Delta(\bar{X}_s, 1-s)\frac{1}{p_{1-s}(\bar{X}_s)}\frac{\mathrm{d}}{\mathrm{d}x}(\sigma^2(1-s,\bar{X}_s)p_{1-s}(\bar{X}_s))\,\mathrm{d}s$$
$$- \int_{1-t}^1 f_\Delta(\bar{X}_s, 1-s)\sigma(1-s,\bar{X}_s)\,\mathrm{d}\bar{W}_s$$

and then

$$E\left(\left|\int_0^t \int_\mathbb{R} f_\Delta(x,s)\,\mathrm{d}L_s^x\right|\right) \leq \|f_\Delta\|_*. \tag{1}$$

Now, given $f \in \mathcal{H}$, let us consider $\{f_n\}_{n \in \mathbb{N}}$, a sequence of elementary functions converging to $f$ in $\mathcal{H}$, and define

$$\int_0^t \int_{\mathbb{R}} f(x,s) \, \mathrm{d}L_s^x = L^1 - \lim_{n \to \infty} \left( \int_0^t \int_{\mathbb{R}} f_n(x,s) \, \mathrm{d}L_s^x \right).$$

Using inequality (1), it is clear that this limit exists and that the definition does not depend on the choice of sequence $(f_n)$.

On the other hand, for any $f \in \mathcal{H}$, it is easy to check that $f(X_\cdot, \cdot)$ is forward and backward integrable with respect to $X$ (see Remark 2.4). Finally, from the uniqueness of the extension, we have

$$\int_0^t \int_{\mathbb{R}} f(x,s) \, \mathrm{d}L_s^x = \int_0^t f(X_s,s) \, \mathrm{d}X_s - \int_0^t f(X_s,s) \, \mathrm{d}^* X_s. \qquad \square$$

**Remark 3.3.** When (H3) is not satisfied, $\|\cdot\|_*$ is only a seminorm (not a norm). In this case, we should restrict our results to $f \in L^2(\mathbb{R} \times [0,1])$. Let us consider the space

$$\mathcal{H} = \{ f \in L^2(\mathbb{R} \times [0,1]) : \|f\|_* < +\infty \}$$

with the norm

$$\|\cdot\| = \|\cdot\|_* + \|\cdot\|_2.$$

The results given in Theorem 3.2 are then true with this new norm and this space.

That is, we are able to define the integral for functions that are not $L^2$, but we then need more regularity over the diffusion.

**Remark 3.4.** Under the hypotheses of Theorem 2.1, if $f \in \mathcal{H}$ then, for any $t \in [0,1]$,

$$\int_0^1 \int_{\mathbb{R}} f(x,s) \, \mathrm{d}L_s^x = -[f(X,\cdot), X]_t.$$

## 4. Itô's formula extension

Our main result states the following:

**Theorem 4.1.** *Assume that $X_t$ is a diffusion process that satisfies statements 1 and 2 of the hypotheses* (H2).

*Let $F$ be a function defined on $\mathbb{R} \times [0,1]$ such that $F$ admits first order Radon–Nikodym derivatives with respect to each parameter and assume that these derivatives satisfy, for every $A \in \mathbb{R}$,*

$$\int_0^1 \int_{-A}^A \left| \frac{\partial F}{\partial t}(x,s) \right| \mathrm{d}x \, r(s) \, \mathrm{d}s < +\infty,$$

$$\int_0^1 \int_{-A}^{A} \left(\frac{\partial F}{\partial x}(x,s)\right)^2 dx\, r(s)\, ds < +\infty.$$

*Then, for all $t \in [0,1]$,*

$$F(X_t, t) = F(X_0, 0) + \int_0^t \frac{\partial F}{\partial x}(X_s, s)\, dX_s + \int_0^t \frac{\partial F}{\partial t}(X_s, s)\, ds - \frac{1}{2}\int_0^t \int_{\mathbb{R}} \frac{\partial F}{\partial x}(x,s)\, dL_s^x.$$

**Proof.** Using localization arguments, we can assume that $\sigma$ and $b$ are bounded, $F$ has compact support and

$$\int_0^1 \int_{\mathbb{R}} \left|\frac{\partial F}{\partial t}(x,s)\right| dx\, r(s)\, ds < +\infty,$$

$$\int_0^1 \int_{\mathbb{R}} \left(\frac{\partial F}{\partial x}(x,s)\right)^2 dx\, r(s)\, ds < +\infty.$$

Also, note that we can suppose, without any loss of generality, that $r(1) \neq 0$.

Let $g \in \mathcal{C}^\infty$ be a function with compact support from $\mathbb{R}$ to $\mathbb{R}^+$ such that $\int_{\mathbb{R}} g(s)\, ds = 1$. We define, for any $n \in \mathbb{N}$, $g_n(s) = n g(ns)$ and

$$F_n(x, t) = \int_0^1 \int_{\mathbb{R}} F(y, s) g_n(t-s) g_n(x-y)\, dy\, ds.$$

Then $F_n \in \mathcal{C}^\infty(\mathbb{R} \times [0,1])$. Hence, by the usual Itô formula, for every $\varepsilon > 0$, we can write

$$F_n(X_t, t) = F_n(X_\varepsilon, \varepsilon) + \int_\varepsilon^t \frac{\partial F_n}{\partial x}(X_s, s)\, dX_s + \int_\varepsilon^t \frac{\partial F_n}{\partial t}(X_s, s)\, ds \quad (2)$$
$$+ \frac{1}{2}\int_\varepsilon^t \sigma^2(s, X_s) \frac{\partial^2 F_n}{\partial x^2}(X_s, s)\, ds.$$

Using the arguments of Azéma *et al.* [1], we will study the convergence of (2).

Since $F$ is a continuous function with compact support, it is easy to check that $(F_n(X_t, t))_{n \in \mathbb{N}}$ converges in probability to $F(X_t, t)$.

On the other hand,

$$\int_0^1 \int_{\mathbb{R}} \left(\frac{\partial F}{\partial x}(x,s)\right)^2 dx\, ds \leq \frac{1}{r(1)} \int_0^1 \int_{\mathbb{R}} \left(\frac{\partial F}{\partial x}(x,s)\right)^2 dx\, r(s)\, ds < +\infty.$$

Hence, $\frac{\partial F}{\partial x} \in L^2(\mathbb{R} \times [0,1])$. Therefore, $(\int_\varepsilon^t \frac{\partial F_n}{\partial x}(X_s, s)\, dX_s)_{n \in \mathbb{N}}$ converges in probability to $(\int_\varepsilon^t \frac{\partial F}{\partial x}(X_s, s)\, dX_s)$. Indeed, using the boundedness of $b$,

$$E\left(\left|\int_\varepsilon^t \left(\frac{\partial F_n}{\partial x}(X_s, s) - \frac{\partial F}{\partial x}(X_s, s)\right) b(s, X_s)\, ds\right|\right)$$

$$\leq K r(\varepsilon) \int_{\varepsilon}^{t} \int_{\mathbb{R}} \left| \frac{\partial F_n}{\partial x}(x,s) - \frac{\partial F}{\partial x}(x,s) \right| dx \, ds \tag{3}$$

and the right-hand side goes to zero when $n$ tends to infinity, since $\frac{\partial F}{\partial x} \in L^2(\mathbb{R} \times [0,1])$ and

$$\frac{\partial F_n}{\partial x}(x,t) = \int_0^1 \int_{\mathbb{R}} \frac{\partial F}{\partial x}(y,s) g_n(t-s) g_n(x-y) \, dy \, ds.$$

Moreover, using the fact that $\sigma$ is a bounded function, we have

$$E\left(\left(\int_{\varepsilon}^{t} \left(\frac{\partial F_n}{\partial x}(X_s,s) - \frac{\partial F}{\partial x}(X_s,s)\right) \sigma(s,X_s) \, dW_s\right)^2\right)$$

$$= E\left(\int_{\varepsilon}^{t} \left(\frac{\partial F_n}{\partial x}(X_s,s) - \frac{\partial F}{\partial x}(X_s,s)\right)^2 \sigma^2(s,X_s) \, ds\right) \tag{4}$$

$$\leq K r(\varepsilon) \int_{\varepsilon}^{t} \int_{\mathbb{R}} \left(\frac{\partial F_n}{\partial x}(x,s) - \frac{\partial F}{\partial x}(x,s)\right)^2 dx \, ds$$

and the final expression converges to zero when $n$ tends to infinity.

Similarly, we obtain that $\frac{\partial F}{\partial t} \in L^1(\mathbb{R} \times [0,1])$ and we prove that $(\int_{\varepsilon}^{t} \frac{\partial F_n}{\partial t}(X_s,s) \, ds)_{n \in \mathbb{N}}$ converges in probability to $(\int_{\varepsilon}^{t} \frac{\partial F}{\partial t}(X_s,s) \, ds)$.

So, letting $n$ tend to infinity in (2), we get that the sequence

$$\left(\frac{1}{2} \int_{\varepsilon}^{t} \sigma^2(s,X_s) \frac{\partial^2 F_n}{\partial x^2}(X_s,s) \, ds\right)_{n \in \mathbb{N}}$$

converges in probability to

$$F(X_t,t) - F(X_\varepsilon,\varepsilon) - \int_{\varepsilon}^{t} \frac{\partial F}{\partial x}(X_s,s) \, dX_s - \int_{\varepsilon}^{t} \frac{\partial F}{\partial t}(X_s,s) \, ds.$$

Note that by using Theorem 2.1 (see also Remarks 2.2 and 2.3) and Theorem 3.2, we can write

$$\int_{\varepsilon}^{t} \sigma^2(s,X_s) \frac{\partial^2 F_n}{\partial x^2}(X_s,s) \, ds$$

$$= \left[\frac{\partial F_n}{\partial x}(X,\cdot),X\right]_t - \left[\frac{\partial F_n}{\partial x}(X,\cdot),X\right]_\varepsilon$$

$$= \int_0^1 \frac{\partial F_n}{\partial x}(X_s,s) I_{(\varepsilon,t)}(s) \, d^*X_s - \int_0^1 \frac{\partial F_n}{\partial x}(X_s,s) I_{(\varepsilon,t)}(s) \, dX_s$$

$$= -\int_0^1 \int_{\mathbb{R}} \frac{\partial F_n}{\partial x}(x,s) I_{(\varepsilon,t)}(s) \, dL_s^x,$$

since $\frac{\partial F_n}{\partial x}(x,s)I_{(\varepsilon,t)}(s) \in \mathcal{H}$.

Moreover, $(\frac{\partial F_n}{\partial x}(x,s)I_{(\varepsilon,t)}(s), x \in \mathbb{R}, s \in [0,1])_{n \in \mathbb{N}}$ converges in $\mathcal{H}$ to $(\frac{\partial F}{\partial x}(x,s)I_{(\varepsilon,t)}(s)$, $x \in \mathbb{R}, s \in [0,1])$. Indeed, we will use the space $\mathcal{H}$ defined in Remark 3.3. Since $\frac{\partial F_n}{\partial x}$ and $\frac{\partial F}{\partial x}$ belong to $L^2(\mathbb{R} \times [0,1])$, it suffices to check the convergence using the seminorm $\|\cdot\|_*$. Now, observe that in (3) and (4), we have proved the convergence to zero of two of the terms appearing in the definition of the seminorm. Hence, it is enough to check that

$$\int_\varepsilon^t \int_\mathbb{R} \left|\frac{\partial F_n}{\partial x}(x,s) - \frac{\partial F}{\partial x}(x,s)\right| \left|\frac{\mathrm{d}}{\mathrm{d}x}(\sigma^2(s,x)p_s(x))\right| \mathrm{d}x\,\mathrm{d}s$$

converges to zero. But this term can be written as

$$E\int_\varepsilon^t \left|\frac{\partial F_n}{\partial x}(X_s,s) - \frac{\partial F}{\partial x}(X_s,s)\right| \frac{1}{p_s(X_s)} \left|\frac{\mathrm{d}}{\mathrm{d}x}(\sigma^2(s,X_s)p_s(X_s))\right| \mathrm{d}s$$

$$\leq \left(E\int_\varepsilon^t \left(\frac{\partial F_n}{\partial x}(X_s,s) - \frac{\partial F}{\partial x}(X_s,s)\right)^2 \mathrm{d}s\right)^{1/2}$$

$$\times \left(E\int_\varepsilon^t \frac{1}{p_s^2(X_s)} \left(\frac{\mathrm{d}}{\mathrm{d}x}(\sigma^2(s,X_s)p_s(X_s))\right)^2 \mathrm{d}s\right)^{1/2}$$

$$\leq (r(\varepsilon))^{1/2}(t-\varepsilon)^{1/2}u(\varepsilon)\left(\int_\mathbb{R}\int_\varepsilon^t \left(\frac{\partial F_n}{\partial x}(x,s) - \frac{\partial F}{\partial x}(x,s)\right)^2 \mathrm{d}s\,\mathrm{d}x\right)^{1/2},$$

the final expression converging to zero when $n$ tends to infinity.

So, we have proved that $(\frac{\partial F_n}{\partial x}(x,s)I_{(\varepsilon,t)}(s), x \in \mathbb{R}, s \in [0,1])_{n \in \mathbb{N}}$ converges in $\mathcal{H}$ to $(\frac{\partial F}{\partial x}(x,s)I_{(\varepsilon,t)}(s), x \in \mathbb{R}, s \in [0,1])$ and that clearly yields that

$$\left(\int_0^1 \int_\mathbb{R} \frac{\partial F_n}{\partial x}(x,s)I_{(\varepsilon,t)}(s)\,\mathrm{d}L_s^x\right)_{n \in \mathbb{N}}$$

converges in $L^1$ to $\int_0^1 \int_\mathbb{R} \frac{\partial F}{\partial x}(x,s)I_{(\varepsilon,t)}(s)\,\mathrm{d}L_s^x$.

Now, we have that for any $\varepsilon > 0$,

$$F(X_t,t) = F(X_\varepsilon,\varepsilon) + \int_\varepsilon^t \frac{\partial F}{\partial x}(X_s,s)\,\mathrm{d}X_s + \int_\varepsilon^t \frac{\partial F}{\partial t}(X_s,s)\,\mathrm{d}s$$
$$- \frac{1}{2}\int_0^1 \int_\mathbb{R} \frac{\partial F}{\partial x}(x,s)I_{(\varepsilon,t)}(s)\,\mathrm{d}L_s^x. \tag{5}$$

Letting $\varepsilon$ tend to zero, the proof is completed. $\square$

**Remark 4.2.** Note that under the hypotheses of Theorem 4.1, it is possible that $\frac{\partial F}{\partial x}$ does not belong to the space $\mathcal{H}$ established in Remark 3.3. In this case, using localization arguments, we can always assume that $(\frac{\partial F}{\partial x}(x,s)I_{(\varepsilon,t)}(s), x \in \mathbb{R}, s \in [0,1])$ belongs to $\mathcal{H}$ for any $\varepsilon > 0$ and can define

$$\int_0^t \int_\mathbb{R} \frac{\partial F}{\partial x}(x,s) \, \mathrm{d}L_s^x = \lim_{\varepsilon \to 0} \int_0^1 \int_\mathbb{R} \frac{\partial F}{\partial x}(x,s) I_{(\varepsilon,t)}(s) \, \mathrm{d}L_s^x.$$

This limit exists in probability since all the other limits in (5) exist.

## 5. Examples

Bardina and Jolis [2] gave conditions on the coefficients which ensure that the hypotheses of Theorem 4.1 over the diffusion are satisfied. For more precise details see Theorems 3.7 and 3.10 of Bardina and Jolis [2] for strongly elliptic diffusions and Theorem 3.9 for elliptic diffusions.

## Acknowledgements

This work was partially supported by MEC-FEDER Grants MTM2006-06427 and MTM2006-01351.